\documentclass{article}
\usepackage{amsmath,psfig}

\begin{document}

\title{New Numerical Algorithm for Modeling of Boson-Fermion Stars in
Dilatonic Gravity\!
\thanks{This work was supported by Bulgarian National
Scientific Fund, Contr. NoNo F610/99, MM602/96 and by Sofia
University Research Fund, Contr. No. 245/99.}}
\author{T.L.Boyadjiev\thanks{Faculty of Mathematics and Computer Science,
University of Sofia, 5 James Bourchier
Blvd., 1164 Sofia, Bulgaria,
E-mail: todorlb@fmi.uni-sofia.bg}\hspace{2mm},
M.D.Todorov\thanks{Faculty of Applied Mathematics
and Computer Science,
Technical University of Sofia, 1756 Sofia, Bulgaria,
E-mail: mtod@vmei.acad.bg}\hspace{2mm},
P.P.Fiziev\thanks{Faculty of Physics, University of
Sofia, 5, James Bourchier Blvd.,
1164 Sofia, Bulgaria,
E-mail: fiziev@phys.uni-sofia.bg}\hspace{2mm},
S.S.Yazadjiev\thanks{Faculty of Physics, University
of Sofia, 5, James Bourchier Blvd.,
1164 Sofia, Bulgaria,
E-mail: yazad@phys.uni-sofia.bg}}
\date{}
\maketitle

\begin{abstract}
We investigate numerically a models of the static spherically symmetric
boson-fermion stars in scalar-tensor theory of gravity with massive
dilaton field. The proper mathematical model of such stars is
interpreted as a nonlinear two-parametric eigenvalue problem with
unknown internal boundary. We employ the Continuous Analogue of Newton
Method (CANM) which leads on each iteration to two separate linear
boundary value problems with different dimensions inside and outside
the star, respectively. Along with them a nonlinear algebraic system
for the spectral parameters - radius of the star $R_{s}$ and quantity
$\Omega $ is solved also.

In this way we obtain the behaviour of the basic geometric quantities
and functions describing dilaton field and matter fields which
build the star.

\medskip
{\bf Keywords.}
two-parametric nonlinear eigenvalue problem, Continuous Analog of the
Newton Method, mixed fermion-boson stars, scalar-tensor theory of
gravity, massive dilaton field.

\medskip
{\bf AMS subject classifications.}
65C20, 65P20, 65P30, 8308, 83D05
\end{abstract}

\section{Introduction}

In this paper we present a numerical method for solving the equations
of the general scalar-tensor theories of gravity including a dilaton
potential term for the general case of mixed boson-fermion star. This
method is an impruvment of the method concerning the similar problem
proposed in our recent work \cite{tomiplas}. The original domain is
splitted to two domains: inner inside the star and external outside the
star. The solutions in these regions are obtained separately and after
that they are matched. Before we go to the substantial part of this
paper we will describe briefly the physical problem.

\section{Main model}\label{sec2}
Boson stars are gravitationally bound macroscopic quantum states
made up of scalar bosons \cite{K,RB,RB1,GW}. They differ from the usual fermionic
stars in that they are only prevented from collapsing
gravitationally by the Heisenberg uncertainty principle. For
self-interacting boson field the mass of the boson star, even for
small values of the coupling constant turns out to be of the order
of Chandrasekhar's mass when the boson mass is similar to a
proton mass. Thus, the boson stars arise as possible candidates
for non-baryonic dark matter in the universe and consequently as a
possible solution of the one of outstanding problems in today's
astrophysics - the problem of nonluminous matter in the universe.
Most of the stars are of primordial origin being formed from an
original gas of fermions and bosons in the early universe. That is
why it should be expected that most stars are a mixture of both,
fermions and bosons in different proportions.

Boson-fermion  stars are also good model to learn more about the
nature of strong gravitational fields not only in general
relativity but also in the other theories of gravity.

The most natural and promising generalizations of general relativity are the
scalar-tensor theories of gravity \cite{BD,D,Will,DEF1}. In these
theories the gravity is mediated not only by a tensor field (the metric of
space-time) but also by a scalar field (the dilaton). These dilatonic
theories of gravity contain arbitrary functions of the scalar field that
determine the gravitational ``constant'' as a dynamical variable and the
strength of the coupling between the scalar field and matter. It should be
stressed that specific scalar-tensor theories of gravity arise naturally as
low energy limit of string theory \cite{GSW,CFMP,FT,SS1,MV,MAHS} which is the
most promising modern model of unification of all fundamental physical
interactions.

Boson stars in scalar-tensor theories of gravity with massless
dilaton have been widely investigated recently both numerically
and analytically \cite{TX,GJ,T,TLS,CS,BS,Y}. Mixed boson-fermion stars in
scalar tensor theories of gravity however have not been
investigated so far in contrast to general relativistic case where
boson-fermion stars have been investigated \cite{HLM}.

In present paper we consider boson-fermion stars in the most
general scalar-tensor theory of gravity with massive dilaton.

In the Einstein frame the field equations in presence of fermion
and boson matter are:
\begin{eqnarray}
G_{m }^{n }=\kappa _{*}\left( \stackrel{B}{T_{m }^{n }}+\stackrel{F}{%
T_{m }^{n }}\right) +2\partial _{m }\varphi \partial
^{n }\varphi &&\!\!\!\!\!\!
-\partial ^{l }\varphi \partial _{l }\varphi \delta _{m }^{n }+%
\frac{1}{2}U(\varphi )\delta _{m }^{n }, \nonumber \\
\nabla _{m }\nabla ^{m }\varphi +{\frac{1}{4}}U^{\prime }(\varphi )&&\!\!\!\!\!\! =-%
\frac{\kappa _{*}}{2}\alpha (\varphi )\left( \stackrel{B}{T}+\stackrel{F}{T}%
\right), \nonumber \\ [-1.8ex]
\label{SFE} \\[-1.8ex]
\nabla _{m }\nabla ^{m }\Psi +2\alpha (\varphi )\partial ^{l
}\varphi \partial _{l }\Psi &&\!\!\!\!\!\! =-2A^{2}(\varphi ){\frac{\partial {\tilde{%
W}}(\Psi ^{+}\Psi )}{\partial \Psi ^{+}}},  \nonumber \\ \nabla
_{m }\nabla ^{m }\Psi ^{+}+2\alpha (\varphi )\partial ^{l
}\varphi \partial _{l }\Psi ^{+}&&\!\!\!\!\!\! =-2A^{2}(\varphi ){\frac{\partial {%
\tilde{W}}(\Psi ^{+}\Psi )}{\partial \Psi ^{+}}}  \nonumber
\end{eqnarray}
where $\nabla _{m }$ is the Levi-Civita connection with respect
to the metric $g_{m n },(m = 0,...,3;$ $n = 0,...,3)$. The
constant $\kappa _{*}$ is given by $\kappa _{*}=8\pi G_{*}$ where
$G_{*}$ is the bare Newtonian
gravitational constant. The physical gravitational ``constant'' is $%
G_{*}A^{2}(\varphi )$ where $A(\varphi )$ is a function of the
dilaton field
$\varphi $ depending on the concrete scalar-tensor theory of gravity. ${\tilde%
W}(\Psi ^{+}\Psi )$ is the potential of boson field.
The
dilaton potential $U(\varphi )$ can be written in the form
$U(\varphi )=m_{D}^{2}V(\varphi )$ where $m_{D}$ is the dilaton
mass and $V(\varphi )$ is dimensionless function of $\varphi $.

The function $\alpha (\varphi )=\frac{d}{d\varphi }\left[ \ln
A(\varphi )\right] $ determines the strength of the coupling
between the dilaton field
$\varphi $ and the matter. The functions $\stackrel{B}{T}$ and $\stackrel{F}{%
T}$ are correspondingly the trace of the energy-momentum
tensor of the fermionic matter\footnote{In the present article
we consider fermionic matter only in macroscopic approximation,
i.e., after averaging quantum fluctuations of the corresponding
fermion fields. Thus we actually consider standard classical
relativistic matter.} $\stackrel{F}{ T_{m }^{n }}$ and
bosonic matter $\stackrel{B}{T_{m }^{n }}$. Their
explicit forms are
\begin{eqnarray}
\stackrel{B}{T_{m }^{n }}=\!\!\!\!\!\!&& \frac{1}{2}A^{2}(\varphi
)\left(
\partial _{m }\Psi ^{+}\partial ^{n }\Psi +\partial _{m
}\Psi \partial ^{n }\Psi ^{+}\right)  \nonumber \\ \!\!\!\!\!\!&& \qquad
-\frac{1}{2}A^{2}(\varphi )\left[ \partial _{l }\Psi ^{+}\partial
^{l }\Psi -2A^{2}(\varphi ){\tilde{W}}(\Psi ^{+}\Psi )\right]
\delta _{m }^{n }\>,  \label{eq:2a} \\ \stackrel{F}{T_{m }^{n
}}=\!\!\!\!\!\!&& \left( \varepsilon +p\right) u_{m }\,u^{n }-p\,\delta _{m
}^{n }\>. \label{eq:2b}
\end{eqnarray}
Here $\Psi $ is a complex scalar field describing the bosonic matter while $%
\Psi ^{+}$ is its complex conjugated function. The energy density
and the pressure of the fermionic fluid in the Einstein frame are
$\varepsilon =A^{4}(\varphi ){\tilde{\varepsilon}}$ and
$p=A^{4}(\varphi )\,\tilde{p}$ where ${\tilde{\varepsilon}}$ and
$\tilde{p}$ are the physical energy density and pressure. Instead
to give the equation of state of the fermionic matter in the form
$\tilde{p}=\tilde{p}({\tilde{\varepsilon}})$ it is more convenient
to write it in a parametric form
\begin{equation}
{\tilde{\varepsilon}}={\tilde{\varepsilon}_{0}}g(\mu )\,\,\,\,\,\,\,\,\,\,%
\tilde{p}={\tilde{\varepsilon}_{0}\,}f(\mu )  \label{eq:2c}
\end{equation}
where ${\tilde{\varepsilon}_{0}}$ is a properly chosen dimensional
constant and $\mu $ is dimensionless Fermi momentum.

The physical four-velocity of the fluid is denoted by $u_{\mu }$.
The potential for the boson field has the form
\[
{\tilde{W}}(\Psi ^{+}\Psi )=-\frac{m_{D}^{2}}{2}\Psi ^{+}\,\Psi -\frac{1}{4%
}{\tilde{\Lambda}\,}(\Psi ^{+}\Psi )^{2}.
\]
Field equations together with the Bianchi identities lead to the
local conservation law of the energy-momentum of matter
\begin{equation}
\nabla _{n }\stackrel{F}{T_{m }^{n }}=\alpha (\varphi )\stackrel{F}{ T}%
\partial _{m }\varphi \,\,\,.  \label{TBIAN}
\end{equation}

We will consider a static and spherically symmetric mixed
boson-fermion star in asymptotic flat space-time. This means that
the metric $g_{m n }$ has the form
\begin{equation}
ds^{2}= \exp \left[ \nu ({\cal R}) \right] dt^{2}- \exp \left[ \lambda
({\cal R}) \right] d{\cal R}^{2}-{\cal R}^{2}\left( d\theta ^{2}+\sin ^{2}\theta
\,d\psi ^{2}\right)  \label{eq:gmn}
\end{equation}
where ${\cal R}, \theta ,\psi $ are usual spherical coordinates. The
field configuration is static when the boson field $\Psi $
satisfies the relationship
\[
\Psi ={\tilde{\sigma}(}{\cal R}{)\,}\exp(i\omega t).
\]
Here $\omega $ is a real number and ${\tilde{\sigma}(}{\cal R})$ is a
real function. Taking into account the assumption that have been
made the system of the field equation is reduced to a system of
ordinary differential equations. Before we explicitly write the
system we are going to introduce a rescaled (dimensionless) radial
coordinate by $r= m_{B}{\cal R}$ where $m_{B}$ is the mass of the
bosons. From now on, a prime will denote a differentiation with
respect to the dimensionless radial coordinate $r$. After introducing the
dimensionless quantities
\[
\Omega ={\frac{\omega }{m_{B}}},\quad \sigma =\sqrt{\kappa }_{*}\,{\tilde{%
\sigma}},\quad \Lambda ={\frac{{\tilde{\Lambda}}}{\kappa _{*}{m^{2}}_{B}}}%
,\quad \gamma ={\frac{m_{D}}{m_{B}}}.
\]
and defining the dimensionless energy-momentum tensors as
$T_{m}^{n}= {\frac {\kappa_{*}} {m_{B}^2}}
T_{m}^{n}$ the components of the dimensionless
energy-momentum tensor of the fermionic and bosonic matter become
correspondingly
\begin{eqnarray*}
&&\!\!\!\!\!\!\qquad\qquad \stackrel{\mathit{F}}{T_{0}^{0}} =bA^{4}(\varphi )\,g(\mu ),\qquad \stackrel{%
\mathit{F}}{T_{1}^{1}}=-bA^{4}(\varphi )\,f(\mu ),\\             
\qquad \stackrel{\mathit{B}}{T_{0}^{0}}&&\!\!\!\!\!\! =\frac{1}{2}A^{2}(\varphi )\left[ \Omega
^{2}\sigma ^{2}(r)\exp (-\nu(r) )+\left( \frac{\partial \sigma }{dr}\right)
^{2}\exp (-\lambda(r) )\,\right] -A^{4}(\varphi )W(\sigma ^{2}),\\ 
\qquad \stackrel{\mathit{B}}{T_{1}^{1}}&&\!\!\!\!\!\! =-\frac{1}{2}A^{2}(\varphi )\left[\Omega
^{2}\sigma ^{2}(r)\exp (-\nu(r) ) + \left( \frac{\partial \sigma }{\partial r}%
\right) ^{2}\exp (-\lambda(r) )\,\right] - A^{4}(\varphi )W(\sigma ^{2}).
\end{eqnarray*}
Here the parameter $b=\frac{\kappa
_{*}{\tilde{\varepsilon}_{0}}}{m_{B}^{2}} $ describes the relation
between the Compton length of dilaton and the usual radius of
neutron star in general relativity.

The functions $\stackrel{\mathit{B}}{T}$ and $\stackrel{\mathit{F}}{T}$
describing the trace of energy-momentum tensor have a form:
\[
\stackrel{\mathit{B}}{T}=-A^{2}(\varphi )\,\left[ \Omega ^{2}\sigma
^{2}(r)\exp (-\nu(r) ) - \left( \frac{\partial \sigma }{dr}\right) ^{2}\exp
(-\lambda(r) )\,\right] - 4A^{4}(\varphi )W(\sigma ^{2}),
\]
\[
\stackrel{\mathit{F}}{T}=bA^{4}(\varphi )\,\left[ g(\mu )-3\,f(\mu )\right]
.
\]

For the independent (dimensionless) radial coordinate $r$ we have
$r\in [0,R_{s}]\cup [R_{s},\infty )$ where $0<R_{s}<\infty $ is
the unknown radius of the fermionic part of the mixed
boson-fermion star.

With all definitions we have given the main system of differential
equations governing the structure of a static and spherically
symmetric boson-fermion star can be written in the following form:

1. In the interior of the fermionic part of the star (the functions in this domain are subscribed
by~$i$)

\begin{eqnarray}
\hspace{1.5cm} \frac{d\lambda }{dr}&&\!\!\!\!\!\! =F_{1,i}\equiv \frac{1-\exp (\lambda )}{r}+r\,\left\{
\exp (\lambda )\left[ \stackrel{F}{T_{0}^{0}}+\stackrel{B}{T_{0}^{0}}+\frac{1%
}{2}\gamma ^{2}V(\varphi )\right] +\left( \frac{d\varphi }{dr}\right)
^{2}\right\} ,  \nonumber \\
\frac{d\nu }{dr}&&\!\!\!\!\!\! =F_{2,i}\equiv -\frac{1-\exp (\lambda )}{r}-r\,\left\{
\exp (\lambda )\left[ \stackrel{F}{T_{1}^{1}}+\stackrel{B}{T_{1}^{1}}+\frac{1}{2}%
\gamma ^{2}V(\varphi )\right]-\left( \frac{d\varphi }{dr}\right) ^{2}\right\} ,
\nonumber \\
\frac{d^{2}\varphi }{dr^{2}}&&\!\!\!\!\!\! =F_{3,i}\equiv -\frac{2}{r}\frac{d\varphi }{dr}%
+\frac{1}{2}\left( F_{1,i}-F_{2,i}\right) \frac{d\varphi }{dr}\nonumber\\ &&\!\!\!\!\!\!\qquad\quad +\frac{1}{2}%
\exp (\lambda )\,\left[ \alpha (\varphi )\,(\stackrel{F}{T}+\stackrel{B}{T})+%
\frac{1}{2}\gamma ^{2}V^{^{\prime }}(\varphi )\right] ,  \label{eqi} \\
\frac{d^{2}\sigma }{dr^{2}}&&\!\!\!\!\!\! =F_{4,i}\equiv -\frac{2}{r}\frac{d\sigma }{dr}%
+\left[ \frac{1}{2}\left( F_{1,i}-F_{2,i}\right) -2\alpha (\varphi )\frac{%
d\varphi }{dr}\right] \frac{d\sigma }{dr}\nonumber\\ &&\!\!\!\!\!\!\qquad\quad
-\sigma \exp (\lambda )\left[ \Omega ^{2}\exp (-\nu )+2A^{2}(\varphi
)W^{^{\prime }}(\sigma ^{2})\right] , \nonumber \\
\frac{d\mu }{dr}&&\!\!\!\!\!\!\! =F_{5,i}\equiv -\frac{g(\mu )+f(\mu )}{f^{^{\prime }}(\mu )%
}\left[ \frac{1}{2}F_{2}+\alpha (\varphi )\frac{d\varphi }{dr}\right] .
\nonumber
\end{eqnarray}

Here $\lambda (r),\nu (r),\varphi (r),\sigma (r)$ and $\mu (r)$ are unknown
functions of $r$ and $\Omega $ is a unknown parameter.
Having in mind the physical assumptions we have to solve
the equations (\ref{eqi}) under following boundary conditions:
\begin{equation}
\lambda (0)=\frac{d\varphi }{dr}(0)=\frac{d\sigma }{dr}(0)=0,\quad \sigma
(0)=\sigma _{c},\quad \mu (0)=\mu _{c}, \label{bcil}
\end{equation}
\begin{equation}
\mu (R_{s})=0  \label{bcir}
\end{equation}

\noindent where $\sigma _{c}$ and $\mu _{c}$ are the values of density of
the bosonic and fermionic matter, respectively at the star's centre. The first three
conditions in (\ref{bcil}) guarantee the nonsingularity of the metrics and
the functions $\lambda (r),$ $\varphi (r),$ $\sigma (r)$ at the star's
centre.

2. In the external domain (subscribed by $e$) there is not fermionic matter,
i.e., one can suppose formally that the function $\mu (r)\equiv 0$ if $x\geq
R_{s}$. The fermionic part of the energy-momentum tensor vanishes
identically also and thus the differential equations with respect to the
rest four unknown functions $\lambda (r),$ $\nu (r),$ $\varphi (r)$ and $\sigma (r)$
are:

\begin{eqnarray}
\hspace{1.5cm}\frac{d\lambda }{dr}&&\!\!\!\!\!\! =F_{1,e}\equiv \frac{1-\exp (\lambda )}{r}+r\,\left\{
\exp (\lambda )\,\left[ \stackrel{B}{T_{0}^{0}}+\frac{1}{2}\gamma
^{2}V(\varphi )\right] +\left( \frac{d\varphi }{dr}\right) ^{2}\right\},
\nonumber \\
\frac{d\nu }{dr}&&\!\!\!\!\!\! =F_{2,e}\equiv -\frac{1-\exp (\lambda )}{r}-r\,\left\{
\exp (\lambda )\,\left[ \stackrel{B}{T_{1}^{1}}+\frac{1}{2}\gamma
^{2}V(\varphi )\right] -\left( \frac{d\varphi }{dr}\right) ^{2}\right\},
\nonumber \\
\frac{d^{2}\varphi }{dr^{2}}&&\!\!\!\!\!\! =F_{3,e}\equiv -\frac{2}{r}\frac{d\varphi }{dr}%
+\frac{1}{2}\left( F_{1,e}-F_{2,e}\right) \frac{d\varphi }{dr}\nonumber\\ &&\!\!\!\!\!\! \qquad +\frac{1}{2}%
\exp (\lambda )\,\left[ \alpha (\varphi )\stackrel{B}{T}+\frac{1}{2}\gamma
^{2}V^{^{\prime }}(\varphi )\right],  \label{eqe} \\
\frac{d^{2}\sigma }{dr^{2}}&&\!\!\!\!\!\! =F_{4,e}\equiv -\frac{2}{r}\frac{d\sigma }{dr}%
+\left[ \frac{1}{2}\left( F_{1,e}-F_{2,e}\right) -2\alpha (\varphi )\frac{%
d\varphi }{dr}\right] \frac{d\sigma }{dr}\nonumber\\             &&\!\!\!\!\!\! \qquad -\sigma \exp (\lambda )\,\left[
\Omega ^{2}\exp (-\nu )+2A^{2}(\varphi )W^{^{\prime }}(\sigma ^{2})\right].
\nonumber
\end{eqnarray}

As it is required by the asymptotic flatness of space-time the boundary conditions at the
infinity are
\begin{equation}
\nu (\infty )=0,\quad \varphi (\infty )=0,\quad \sigma (\infty )=0
\label{bcer}
\end{equation}
where we denote $(\cdot) (\infty) = \lim_{r \to \infty} (\cdot) (r)$.

We seek for a solution $\left[ \lambda (r),\nu (r),\varphi (r),\sigma
(r),\mu (r),R_{s},\Omega \right] \ $\ subjected to the nonlinear ODEs (\ref
{eqi})\ and (\ref{eqe}),\ satisfying the boundary conditions (\ref{bcil}), (%
\ref{bcir})\ and (\ref{bcer}). At that we assume the function $\mu (r)$ is
continuous in the interval $[0,R_{s}],$ whilst the functions $\lambda (r),$ $%
\nu (r)$\ are continuous and the functions $\varphi (r),$ $\sigma (r)$ are
smooth in the whole interval $[0,\infty )$, including the unknown inner
boundary $r=R_{s}$.

The so posed BVP is a two-parametric eigenvalue problem with respect to the
quantities $R_{s}$ and $\Omega $.

Let us emphasize that a number of methods for solving the free-boundary
problems are considered in detail in \cite{Vab,numrec}.

Here we aim at applying the new solving method to the above formulated
problem. This method differs from that one proposed in \cite{tomiplas} and
for the governing field equations written in the forms (\ref{eqi})\ and (\ref{eqe}%
)\ it possesses certain advantage.

\section{Method of solution}
At first we scale the variable
$r$ using the Landau transformation \cite{Vab} and in this way we obtain a fixed computational domain. Namely
\[
x=\frac{r}{R_{s}},x\in [0,1]\cup [1,\infty ).
\]
For our further considerations it is convenient to present the systems (\ref
{eqi}) and (\ref{eqe})\ in following equivalent form as systems of first
order ODEs:
\begin{eqnarray}
-\mathbf{y}_{i}^{\prime }+R_{s}\mathbf{F}_{i}(R_{s}x,\mathbf{y}_{i},\Omega )
&=&0,  \label{odei} \\
-\mathbf{y}_{e}^{\prime }+R_{s}\mathbf{F}_{e}(R_{s}x,\mathbf{y}_{e},\Omega )
&=&0             \label{odee}
\end{eqnarray}

\noindent with respect to the unknown vector functions $$\mathbf{y}%
_{i}(x)\equiv (\lambda (x),\nu (x),\varphi (x),\xi (x),\sigma (x),\eta
(x),\mu (x))^{T},$$ $$\mathbf{y}_{e}(x)\equiv (\lambda (x),\nu (x),\varphi
(x),\xi (x),\sigma (x),\eta (x))^{T}$$ and right hand sides $\mathbf{F}%
_{i}\equiv (F_{1},F_{2},\xi ,F_{3},\eta ,F_{4},F_{5})^{T}$, $\mathbf{F}%
_{e}\equiv (F_{1},F_{2},\xi ,F_{3},\eta ,F_{4})^{T}$ where $(.)^{\prime }$
stands for differentiation towards the new variable $x$.

For given values of the parameters $R_{s}$\ and $\Omega $\ the independent
solving of the inner system (\ref{odei}) requires seven boundary conditions. At
the same time we have at disposal only six conditions of the kind (\ref{bcil})
and (\ref{bcir}). In order to complete the problem we set additionally one
more parametric condition (the value of someone from among the functions $%
\lambda (x),\nu (x),\varphi (x),\xi (x),\sigma (x)$ or $\eta (x)$) at the
point $x=1)$. Let us set for example
\begin{equation}
\varphi _{i}(1)=\varphi _{s}             \label{dopusl}
\end{equation}
where $\varphi _{s}$\ is a parameter. Then the boundary conditions (\ref
{bcil}), (\ref{bcir}) and (\ref{dopusl}) of the inner BVP can be presented
in the form:

\begin{equation}
B_{0,i}\mathbf{y}_{i}(0)-D_{0,i}=0,\quad B_{1,i}\mathbf{y}%
_{i}(1)-D_{1,i}(\varphi _{s})=0.  \label{bci}
\end{equation}
Here the matrices $B_{0,i}=diag(1,0,0,1,1,1,1)$, $D_{0,i}=diag(0,0,0,0,%
\sigma _{c},0,\mu _{c})$, $B_{1,i}=diag(0,0,1,0,0,0,1)$, $%
D_{1,i}=diag(0,0,\varphi _{s},0,0,0,0)$.

Obviously the solution in the inner domain $x\in [0,1]$ depends not only on
the variable $x,$ but it is a function of the three parameters $R_{s},\Omega
,\varphi _{s}$ as well, i.e., $\mathbf{y}_{i}=\mathbf{y}_{i}(x,\Omega
,R_{s},\varphi _{s}).$

In the external domain $x\geq 1$ the vector of solutions $$\mathbf{y}%
_{e}(x)\equiv (\lambda (x),\nu (x),\varphi (x),\xi (x),\sigma (x),\eta
(x))^{T}$$ is 6D. Thereupon six boundary conditions are indispensable to
solving of the equation (\ref{odee}). At the same time the three boundary
conditions (\ref{bcer}) are only known. Let us consider that the solution $%
\mathbf{y}_{i}(x)$ in the inner domain $x\in [0,1]$ is knowledged. Then
we postulate the rest three deficient conditions to be the continuity
conditions at the point $x=1.$The first of them is similar to the
condition (\ref
{dopusl}) and the else two we assign to arbitrary two functions from among $%
\lambda (x),$ $\nu (x),\mathbf{\ }\xi (x),$ $\sigma (x)$ and $\eta (x)$, for
example
\[
\lambda _{e}(1)=\lambda _{i}(1),\quad \varphi _{e}(1)=\varphi _{s},\quad
\sigma _{e}(1)=\sigma _{i}(1).
\]
It is convenient to present the boundary conditions in the external domain
in matrix form again:
\begin{equation}
B_{1,e}\mathbf{y}_{e}(1)-D_{1,e}(\varphi _{s})=0,\quad B_{\infty ,e}\mathbf{y%
}_{e}(\infty )=0               \label{bce}
\end{equation}
where the matrices $B_{1,e}=diag(1,0,1,0,1,0)$, $D_{1,e}=diag(\lambda
_{i}(1),0,\varphi _{s},0,\sigma _{i}(1),0)$, $B_{\infty
,e}=diag(0,1,1,0,1,0) $.

\medskip Let the solutions $\mathbf{y}_{i}=\mathbf{y}_{i}(x,\Omega
,R_{s},\varphi _{s})$ and $\mathbf{y}_{e}=\mathbf{y}_{e}(x,\Omega
,R_{s},\varphi _{s})$ be supposed known. Generally speaking for given
arbitrary values of the parameters $R_{s},\Omega $ and $\varphi _{s}$ the
continuity conditions with respect to the functions $\nu (x),\xi (x)$ and $%
\eta (x)$ at the point $x=1$ are not satisfied. We choose the parameters $%
R_{s},\Omega $ and $\varphi _{s}$ in such manner that the continuity
conditions for the functions $\nu (x)$, $\xi (x)$ and $\eta (x)$ to be held,
i.e.,

\begin{eqnarray}
\nu _{e}(1,R_{s},\Omega ,\varphi _{s})-\nu _{i}(1,R_{s},\Omega ,\varphi
_{s})&&\!\!\!\!\!\! =0,  \nonumber \\
\xi _{e}(1,R_{s},\Omega ,\varphi _{s})-\xi _{i}(1,R_{s},\Omega ,\varphi
_{s})&&\!\!\!\!\!\! =0,  \label{cc} \\
\eta _{e}(1,R_{s},\Omega ,\varphi _{s})-\eta _{i}(1,R_{s},\Omega ,\varphi
_{s})&&\!\!\!\!\!\! =0.  \nonumber
\end{eqnarray}

These conditions should be interpreted as three nonlinear algebraic
equations in regard to the unknown quantities $R_{s},\Omega $ and $\varphi
_{s}$. The usual way for solving of the above mentioned kind of equations (%
\ref{cc})\ is by means of various iteration methods, for example Newton's
methods. The traditional technology similarly to the methods like shutting
\cite{numrec1}, requires separate treatment of the BVPs and the algebraic
continuity equations and brings itself to additional linear ODEs for
elements of the corresponding to (\ref{cc}) Jacobi matrix. These elements
are functions of the variable $\ x$ and they have to be known actually only
at the point $x=1$. The solving of the both nonlinear BVPs (\ref{odei}), (%
\ref{bci})\ and (\ref{odee}), (\ref{bce}) along with the attached linear
equations is another hard enough task.

At the present work using the CANM \cite{gavurin}, \cite{jmp} we propose a
common treatment of both, differential and algebraic problems.

\mathstrut We suppose that the nonlinear spectral problem (\ref{odei}),
(\ref {bci}), (\ref{odee}), (\ref{bce}) and (\ref{cc}) has a ``well
separated''
\cite{jmp} exact solution. Let the functions $\mathbf{y}_{i,0}(x),\mathbf{y%
}_{e,0}(x)$ and the parameters $R_{s,0},\Omega _{0},\varphi _{s,0}$ are
initial approximations to this solution. The CANM leads to the following
iteration process:
\begin{eqnarray}
\mathbf{y}_{i,k+1}(x) &&\!\!\!\!\!\!=\mathbf{y}_{i,k}(x)+\tau _{k}\mathbf{z}_{i,k}(x),
\label{stepi} \\
\mathbf{y}_{e,k+1}(x) &&\!\!\!\!\!\!=\mathbf{y}_{e,k}(x)+\tau _{k}\mathbf{z}_{e,k}(x),
\label{stepe}\\
R_{s,k+1} &&\!\!\!\!\!\!=R_{s,k}+\tau _{k}\rho _{k},               \label{stepro} \\
\Omega _{k+1}&&\!\!\!\!\!\!=\Omega _{k}+\tau _{k}\omega _{k}, \label{stepomega} \\
\varphi _{s,k+1}&&\!\!\!\!\!\!= \varphi _{s,k}+\tau _{k}\phi _{k}.  \label{stepfis}
\end{eqnarray}

Here $\tau _{k}\in (0,1]$\ is a parameter which can rule the
convergence of iteration process. The increments $\mathbf{z}_{i,k}(x),$ $\mathbf{z}%
_{e,k}(x),\rho _{k},\omega _{k}$ and $\phi _{k}$, $k=0,1,2,...$ satisfy the
linear ODEs (for sake of simplicity henceforth we will omit the number of
iterations $k$):
\begin{eqnarray}
-\mathbf{z}_{i}^{\prime }+R_{s}\frac{\partial \mathbf{F}_{i}}{\partial
\mathbf{y}_{i}}\mathbf{z}_{i}+\left( R_{s}\frac{\partial \mathbf{F}_{i}}{%
\partial R_{s}} + \mathbf{F}_{i}\right) \rho +R_{s}\frac{\partial \mathbf{F}%
_{i}}{\partial \Omega }\omega &&\!\!\!\!\!\!=\mathbf{y}_{i}^{\prime
}-R_{s}\mathbf{F}_{i}, \label{canmi}
\\-\mathbf{z}_{e}^{\prime }+R_{s}\frac{\partial
\mathbf{F}_{e}}{\partial
\mathbf{y}_{e}}\mathbf{z}_{e}+\left( R_{s}\frac{\partial \mathbf{F}_{e}}{%
\partial R_{s}}+\mathbf{F}_{e}\right) \rho +R_{s}\frac{\partial \mathbf{F}%
_{e}}{\partial \Omega }\omega &&\!\!\!\!\!\!=\mathbf{y}_{e}^{\prime
}-R_{s}\mathbf{F}_{e}. \label{canme}
\end{eqnarray}

\noindent All the coefficients and right hand sides as well in the
above two equations are known functions of the arguments $x$, $R_s$,
$\Omega$ by means of the solution from the previous iteration.
We seek for the unknowns $\mathbf{z}_{i}(x)$ of equation (\ref{canmi}) and $%
\mathbf{z}_{e}(x)$ of equation (\ref{canme}) as a linear combinations
with coefficients $\rho ,\omega $ and $\phi $:
\begin{eqnarray}
\mathbf{z}_{i}(x) &=&\mathbf{s}_{i}(x)+\rho \,\mathbf{u}_{i}(x)+\omega
\mathbf{v}_{i}(x)+\phi \mathbf{w}_{i}(x),  \label{decompi} \\
\mathbf{z}_{e}(x) &=&\mathbf{s}_{e}(x)+\rho \,\mathbf{u}_{e}(x)+\omega
\mathbf{v}_{e}(x)+\phi \mathbf{w}_{e}(x).  \label{decompe}
\end{eqnarray}
Here $\mathbf{s}_{i}(x),\mathbf{u}_{i}(x)$, $\mathbf{v}_{i}(x)$, $%
\mathbf{w}_{i}(x)$, $\mathbf{s}_{i}(x)$ and $\mathbf{s}_{e}(x),\mathbf{u}%
_{e}(x)$, $\mathbf{v}_{e}(x)$, $\mathbf{w}_{e}(x)$\ are new unknown
functions, which are defined in either, internal or external domains.
Substituting for the decomposition (\ref{decompi}) into equation (\ref{canmi}%
) after reduction we obtain:
\begin{eqnarray}
-\mathbf{s}_{i}^{\prime
}+Q_{i}(x)\,\mathbf{s}_{i}&&\!\!\!\!\!\!=\mathbf{y}_{i}^{\prime}-R_{s}\mathbf{F}_{i},
\nonumber \\
-\mathbf{u}_{i}^{\prime }+Q_{i}(x)\,\mathbf{u}_{i}&&\!\!\!\!\!\!=-\left( \mathbf{F}%
_{i}+R_{s}\frac{\partial \mathbf{F}_{i}}{\partial R_{s}}\right) ,
\nonumber \\[-1.8ex]
\label{lini}\\[-1.8ex]
-\mathbf{v}_{i}^{\prime }+Q_{i}(x)\,\mathbf{v}_{i}&&\!\!\!\!\!\!=-R_{s}\frac{\partial
\mathbf{F}_{i}}{\partial \Omega },           \nonumber \\
-\mathbf{w}_{i}^{\prime }+Q_{i}(x)\,\mathbf{w}_{i}&&\!\!\!\!\!\!=0  \nonumber
\end{eqnarray}
where $Q_{i}\left( x\right) \equiv R_{s}\frac{\partial \mathbf{F}_{i}\left(
R_{s}x,\mathbf{y}_{i},\Omega \right) }{\partial \mathbf{y}_{i}}$ stands for
a square matrix ($7\times 7)$, which consists of the Frechet derivatives of
operator $\mathbf{F}_{i}$ at the point $\left\{ \mathbf{y}%
_{i}(x),R_{s},\Omega \right\} $.

Similarly applying the CANM to the boundary conditions (\ref{bci}) and
taking into account the dependence of matrix $D_{1,i}$ on the parameter $%
\varphi _{s}$ yields:
\[
B_{0,i}\mathbf{z}_{i}(0)=D_{0,i}-B_{0,i}\mathbf{y}_{i}(0),\quad B_{1,i}%
\mathbf{z}_{i}(1)=D_{1,i}-B_{1,i}\mathbf{y}_{i}(1)-D_{1,i}^{\prime }\phi .
\]
By means of the decomposition (\ref{decompi}) we obtain the following eight
boundary conditions (four left + four right) for the equations (\ref{lini}):
\begin{eqnarray}
&&\!\!\!\!\!\!B_{0,i}\mathbf{s}_{i}(0)=D_{0,i}-B_{0,i}\mathbf{y}_{i}(0),\qquad B_{1,i}%
\mathbf{s}_{i}(1)=D_{1,i}-B_{1,i}\mathbf{y}_{i}(1),\nonumber \\
&&\!\!\!\!\!\!B_{0,i}\mathbf{u}_{i}(0)=0,\qquad\qquad\qquad\qquad\> B_{1,i}\mathbf{u}_{i}(1)=0,\nonumber \\[-1.8ex]
\label{lbci}\\[-1.8ex]
&&\!\!\!\!\!\!B_{0,i}\mathbf{v}_{i}(0)=0,\qquad\qquad\qquad\qquad\> B_{1,i}\mathbf{v}_{i}(1)=0,\nonumber \\
&&\!\!\!\!\!\!B_{0,i}\mathbf{w}_{i}(0)=0,\qquad\qquad\qquad\qquad B_{1,i}\mathbf{w}_{i}(1)=-D_{1,i}^{\prime
}(\varphi _{s})\nonumber .
\end{eqnarray}

Let us now substitute for decomposition (\ref{decompe}) into the linear
equations for external domain (\ref{canme}). As result we obtain the
following four vector equations with regard to the unknown functions $\mathbf{s}%
_{e}(x),\mathbf{u}_{e}(x)$, $\mathbf{v}_{e}(x)$ and $\mathbf{w}_{e}(x)$ with
eight boundary conditions (four left + four right):
\begin{eqnarray}
-\mathbf{s}_{e}^{\prime }+Q_{e}(x)\mathbf{s}_{e}&&\!\!\!\!\!\!=\mathbf{y}_{e}^{\prime
}-R_{s}\mathbf{F}_{e},       \nonumber \\
-\mathbf{u}_{e}^{\prime }+Q_{e}(x)\,\mathbf{u}_{e}&&\!\!\!\!\!\!=-\left( \mathbf{F}%
_{e}+R_{s}\frac{\partial \mathbf{F}_{e}}{\partial R_{s}}\right) ,\nonumber\\[-1.8ex]
\label{line} \\[-1.8ex]
-\mathbf{v}_{e}^{\prime }+Q_{e}(x)\mathbf{v}_{e}&&\!\!\!\!\!\!=-R_{s}\frac{\partial
\mathbf{F}_{e}}{\partial \Omega },           \nonumber \\
-\mathbf{w}_{e}^{\prime }+Q_{e}(x)\,\mathbf{w}_{e}&&\!\!\!\!\!\!=0.  \nonumber
\end{eqnarray}
Here $Q_{e}\left( x\right) \equiv R_{s}\frac{\partial \mathbf{F}%
_{e}\left[ R_{s}x,\mathbf{y}_{e}(x),\Omega \right] }{\partial \mathbf{y}_{e}}
$ is a square matrix $(6\times 6)$ whose elements are Frechet's derivatives
of the operator $\mathbf{F}_{e}$ at the point $\left\{ \mathbf{y}%
_{e}(x),R_{s},\Omega \right\} $.

The corresponding linear BC are obtained in the same way as (\ref{lbci})\
and they become:
\begin{eqnarray}
B_{1,e}\mathbf{s}_{e}(1)&&\!\!\!\!\!\!=D_{1,e}-B_{1,e}\mathbf{y}_{e}(1),\qquad\> B_{\infty ,e}%
\mathbf{s}_{e}(\infty )=-B_{\infty ,e}\mathbf{y}_{e}(\infty ),\nonumber\\
B_{1,e}\mathbf{u}_{e}(1)&&\!\!\!\!\!\!=0,\qquad\qquad\qquad\qquad\quad B_{\infty ,e}\mathbf{u}_{e}(\infty )=0,\nonumber\\[-1.8ex]
\label{lbce}\\[-1.8ex]
B_{1,e}\mathbf{v}_{e}(1)&&\!\!\!\!\!\!=0,\qquad\qquad\qquad\qquad\quad B_{\infty ,e}\mathbf{v}_{e}(\infty )=0,\nonumber\\
B_{1,e}\mathbf{w}_{e}(1)&&\!\!\!\!\!\!=-D_{1,e}^{\prime }(\varphi _{s}),\qquad\qquad\quad B_{\infty ,e}%
\mathbf{w}_{e}(\infty )=0.\nonumber
\end{eqnarray}

In the end to compute the increments $\left\{ \rho ,\omega ,\phi \right\} $
of parameters $R_{s},\Omega $ and $\varphi _{s}$ we use the three conditions
(\ref{cc}).

Let the solutions of linear BVP ($\ref{lini}$), ($\ref{lbci}$) and ($\ref
{line}$), ($\ref{lbce}$) at the $k$th iteration stage are assumed to be
known. For sake of the simplicity we introduce the vector ${\mathbf{\tilde{y}%
(}}x{\mathbf{)}}\equiv (\nu(x) ,\xi(x) ,\eta(x) )^{T}$. For two
arbitrary functions
$h_{i}(x)$ and $h_{e}(x)$, defined in left and right vicinity of the point $x=1$, we set $%
\Delta h\equiv h_{e}(1)-h_{i}(1)$. Then applying the CANM to the equations (%
\ref{cc}) and having in mind the decompositions (\ref{decompi}),(\ref
{decompe}), we attain the vector equation
\begin{equation}
\Delta \mathbf{\tilde{u}\,}\rho +\Delta \mathbf{\tilde{v}\,}\omega +\Delta
\mathbf{\tilde{w}}\,\phi =-\left( \Delta \mathbf{\tilde{y}}+\Delta \mathbf{%
\tilde{s}}\right) ,  \label{lcc}
\end{equation}
which represents an algebraic system consisting of three linear scalar
equations with respect to the three unknowns $\rho $, $\omega $ and
$\phi $.

The general sequence of the algorithm can be recapitulated in the following
way. Let us assume that the functions $\mathbf{y}_{i,k}(x)$, $\mathbf{y}%
_{e,k}(x)$, and parameters $R_{s,k}$, $\Omega _{k}$, $\varphi _{s,k}$ are
given for $k\geq 0$. We solve the linear BVPs (\ref{lini}), (\ref{lbci}) and
thus we compute the functions $\mathbf{s}_{i,k}(x),\mathbf{u}_{i,k}(x)$, $%
\mathbf{v}_{i,k}(x)$, $\mathbf{w}_{i,k}(x)$ in the inner domain $x\in [0,1]$%
. Then we solve the linear BVPs (\ref{line}), (\ref{lbce}) in the external
domain $x\in [1,\infty ]$ and compute the functions $\mathbf{s}_{e,k}(x),%
\mathbf{u}_{e,k}(x)$, $\mathbf{v}_{e,k}(x)$ and $\mathbf{w}_{e,k}(x)$.
Next, to obtain the increments ${\rho }_{k}$, ${\omega }_{k}$ and
${\phi }_{k}$ we solve the linear algebraic system (\ref{lcc}). So
using the decompositions (\ref {decompi}), (\ref{decompe}) and then the
formulae (\ref{stepi}) - (\ref
{stepfis}) we calculate the functions $\mathbf{y}_{i,k+1}(x)$, $\mathbf{y}%
_{e,k+1}(x)$, the radius of the star $R_{s,k+1}$, the quantity $\Omega _{k+1}
$ and the parameter boundary condition $\varphi _{s,k+1}$ as well at the new
iteration stage $k+1$.

At the every iteration $k$\ an optimal time step $\tau _{opt}$ is
determinated in accordance with the Ermakov\&Kalitkin formula
\cite{ermakov}
\[
\tau _{opt}\approx \frac{\delta (0)}{\delta (0)+\delta (1)}
\]
where the residual $\delta (\tau )$ is calculated as follows
\[
\delta (\tau _{k})=\max \,\left[ \delta _{f},(R_{s,k}+\tau _{k}\rho
_{k})^{2},(\Omega _{k}+\tau _{k}\omega _{k})^{2},(\varphi _{s,k}+\tau
_{k}\phi _{k})^{2}\right]
\]
and $\delta _{f}$ is the Euclidean residual of right hand side of the first
equations in the systems (\ref{lini}), (\ref{lbci}) and (\ref{line}), (\ref
{lbce}).

The criterion for termination of the iterations is $\delta (\tau
_{opt})<\varepsilon $ where $\varepsilon \sim 10^{-8}\div 10^{-12}$ for
some $k$.

Taking into account the smoothness of sought solutions we solve the linear BVPs (\ref{lini}), (\ref{lbci}) and (\ref{line}%
), (\ref{lbce})
employing Hermitean splines and spline collocation scheme of
fourth order of approximation \cite{zavyal}.
At that
we utilize essentially the important feature that everyone of the above
mentioned two groups vector BVPs (inner and external) have one and the
same left hand sides.

It is worth to note that the algebraic systems of linear equations
and the system (\ref{lcc})\ as well become ill-posed in the
vicinity of the ``exact'' solution, i.e., for sufficiently small residuals $%
\delta $. That is why for small $\delta $, for example if $\delta <10^{-3}$
(then $\tau _{opt}\sim 1$ usually), it is expedient to use the
Newton-Kantorovich method when the respective matrices are fixed for some $%
\delta \geq 10^{-3}$.

\section{Some numerical results}
\medskip For a purpose of illustrating we will consider and discuss some
results obtained from numerical experiments. A detailed description and
analysing of results from physical point of view will be object of another
our paper.

In present article we consider concrete scalar-tensor model with functions (see Section \ref{sec2})
\begin{eqnarray*}
\qquad\qquad &&\!\!\!\!\!\! A(\varphi )=\exp (\frac{\varphi }{\sqrt{3}}),\qquad V(\varphi )={\frac
{3} {2}}(1-A^{2}(\varphi ))^2,\\ 
f(\mu )&&\!\!\!\!\!\! =\frac{1}{8}\left[ (2\mu -3)\sqrt{\mu +\mu ^{2}}+3\ln \left( \sqrt{%
\mu }+\sqrt{1+\mu }\right) \right],\\  
\;g(\mu )&&\!\!\!\!\!\! =\frac{1}{8}\left[ (6\mu +3)\sqrt{\mu +\mu ^{2}}-3\ln \left(
\sqrt{\mu }+\sqrt{1+\mu }\right) \right],  
\end{eqnarray*}
\[
W(\sigma ^{2})=-\frac{1}{2}\left( \sigma ^{2}+\frac{1}{2}\Lambda \sigma
^{4}\right).
\]
The quantities $b,\Lambda $ are given parameters.
For completeness
we note that in the concrete case the functions $f(\mu)$ and
$g(\mu)$ represent in parametric form  the equation of state of
noninteracting neutron gas while the function $W(\sigma^2)$ describes the boson
field with quartic self-interaction.

\begin{figure}[ht]
\centerline{\psfig{figure=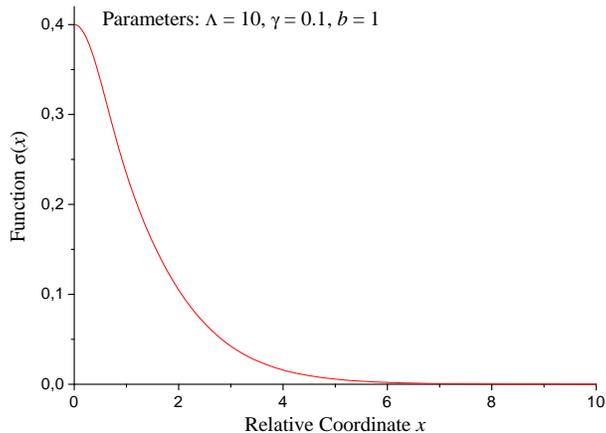,width=4.in}} \caption{The
function $\sigma(x)$ for $\sigma_c = 0.4$; $\mu_c =1.2$.}
\label{sigma}
\end{figure}

\begin{figure}
\centerline{\psfig{figure=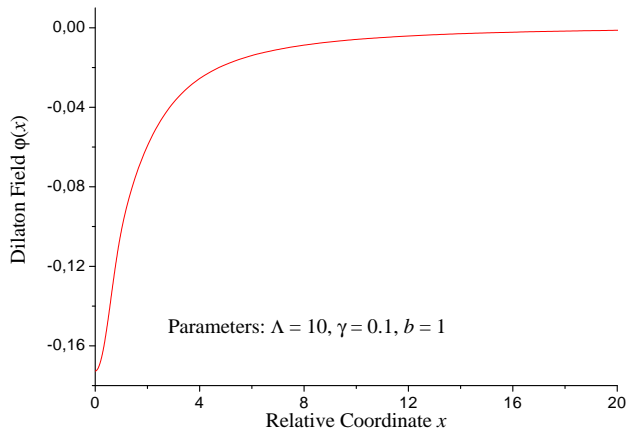,width=4.in}} \caption{The
function $\varphi(x)$ for $\sigma_c = 0.4$; $\mu_c = 1.2$.}
\label{phi}
\end{figure}

\begin{figure}
\centerline{\psfig{figure=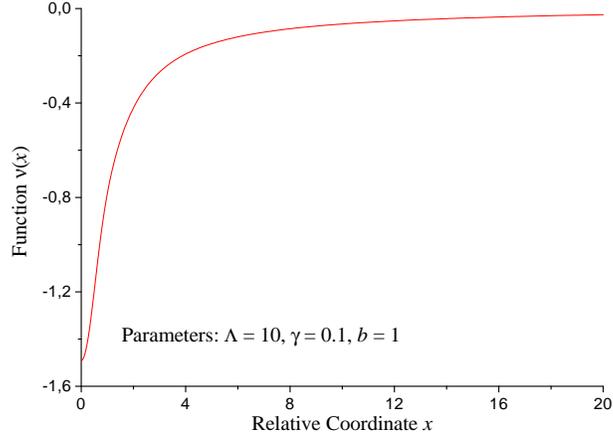,width=4.in}} \caption{The
function $\nu(x)$ for $\sigma_c = 0.4$; $\mu_c = 1.2$.} \label{nu}
\end{figure}

\begin{figure}[ht]
\centerline{\psfig{figure=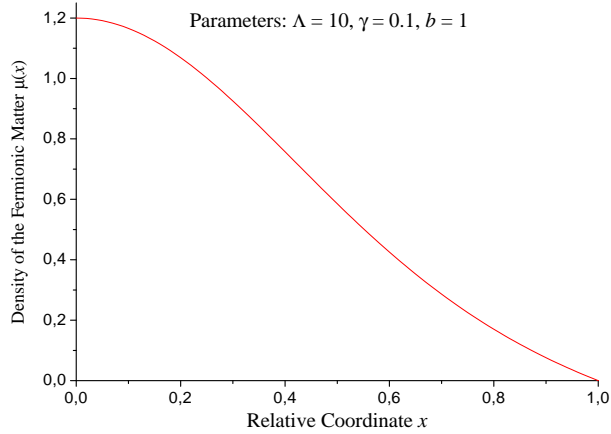,width=4.in}} \caption{The
function $\mu(x)$ for $\sigma_c = 0.4$; $\mu_c = 1.2$.} \label{mu}
\end{figure}

The calculated functions $\sigma(x)$, $\varphi(x)$, $\nu(x)$ and $\mu(x)$ are plotted
correspondingly in Fig. \ref{sigma}, \ref{phi}, \ref{nu} and \ref{mu} for the values of the
parameters $\gamma = 0.1$, $\Lambda =10$ and $b=1.$ The behaviour of the
mentioned functions is typical for wider range of the parameters not only
for those values presented in the figures. The function $\sigma(x)$ decreases
rapidly from its central value $\sigma_c = 0.4$ (in the case under
consideration) to zero, at that when dimensionless coordinate $x > 6$ the
function does not exceed $10^{-4}$. Similarly the function $\nu(x)$ has most
large derivative for $x \in (0,9)$ after that it approaches slowly
zero at infinity like $\frac{1}{x}$. For example when $x \approx 9$ the
derivative
$\nu^{\prime}(x) \approx 10^{-2}$, while for $x > 27$  $\nu^{\prime}(x) < 10^{-4}$,
i.e., the asymptotical behaviour of calculated grid function and its derivative agrees very
well with the theoretical prediction (see \cite{tomiplas}).
The function $\varphi(x)$ increases rapidly for $x<4$, besides that it trends
asymptotically to zero. Obviously the quantitative behaviour of $\varphi(x)$
for central value
$\sigma_c = 0.4$ is determinated by the dominance
of the term $\stackrel{B}{T}$ over the term
$\stackrel{F}{T}$ (see \cite{tomiplas}). At last the function
$\mu(x)$ is nontrivial in the inner domain $x \in [0,1]$, i.e. inside the star.
Here it varies monotonously and continuously from its central value (in the
case under consideration) $\mu_c = 1.2$ till zero at $x = 1$ corresponding to
the radius of the star.

From physical point of view it is important to know the mass of the
boson-fermion star and the total number of particles (bosons and fermions)
making up the star.

The dimensionless star mass can be calculated via the formula
$$
M= \int_{0}^{\infty} r^2 \left(\stackrel{\mathit{B}}{T_{0}^{0}}+
\stackrel{\mathit{F}}{T_{0}^{0}} + \exp(-\lambda)\left({\frac{d\varphi }{dr}}%
\right)^2 + {\frac{\gamma^2}{2}}V(\varphi) \right) dr.
$$

The dimensionless rest mass of the bosons (total number of bosons times the
boson mass) is given by
$$
M_{RB}= \Omega \int_{0}^{\infty} r^2 A^2(\varphi) \exp\left({\frac{\lambda - \nu }{%
2}}\right)\sigma^2  dr.
$$

The dimensionless rest mass of the fermions is correspondingly
\[
M_{RF}=b\int_{0}^{\infty }r^{2}A^{3}(\varphi ) \exp\left({\frac{\lambda }{2}}\right) n(\mu ) dr
\]
where $n(\mu )$ is the density of the fermions. In the case we consider we
have $n(\mu )=\mu ^{\frac{3}{2}}(x)$.

\begin{figure}[ht]
\centerline{\psfig{figure=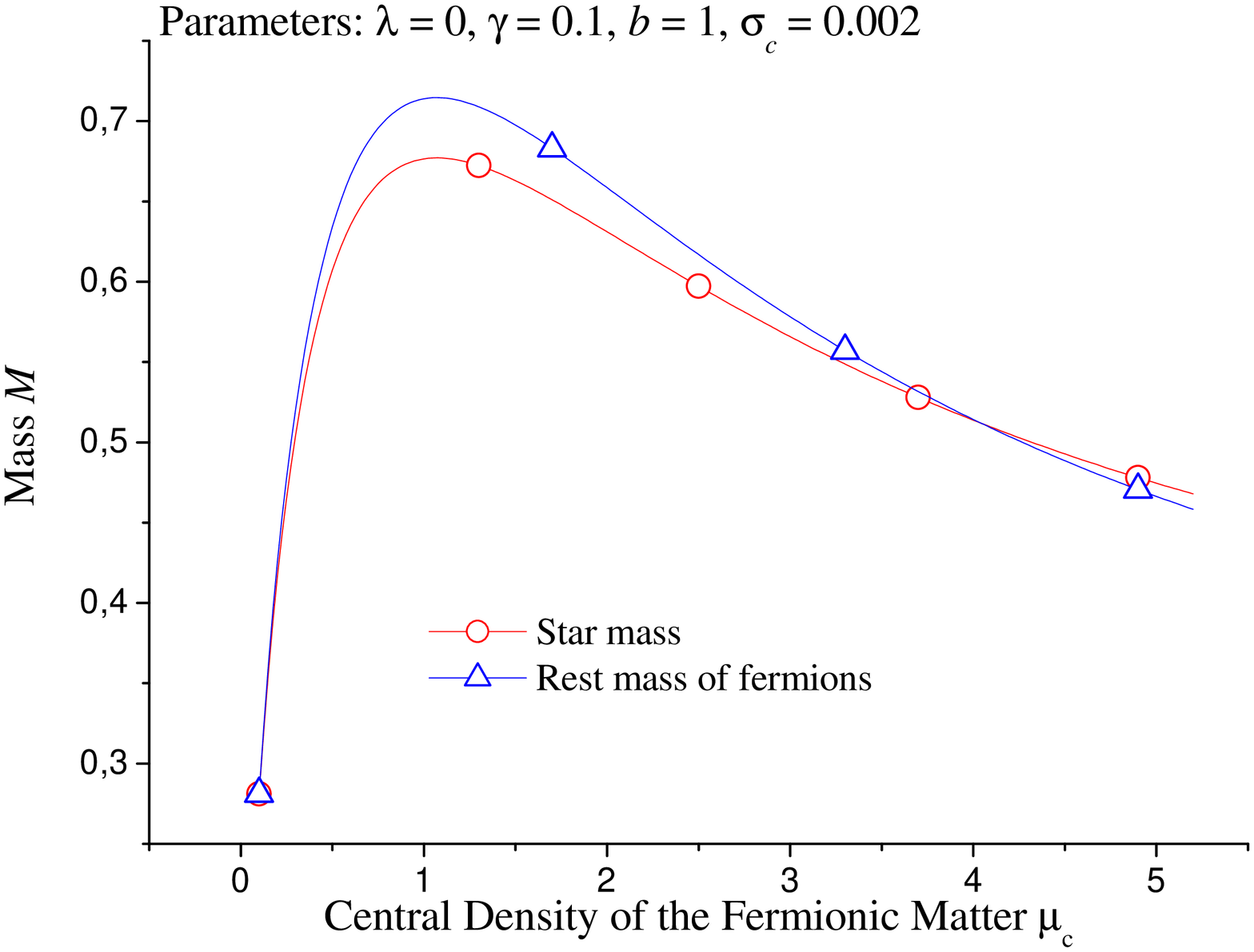,width=4.in}} \caption{The star
mass $M$ and the rest fermion mass $M_{RF}$ as functions of the
central value $\mu_c$} \label{mass}
\end{figure}

The dependencies of the star mass $M$ (solid line) and the rest mass of
fermions $M_{RF}$ (dash line) on the central value $\mu _{c}$ of the
function $\mu(x)$ are shown in configuration diagram on Fig. \ref{mass} for $\lambda =0$,
$\gamma =0.1$, $b=1$ and $\sigma _{c}=0.002$.
It should be
pointed that for so small central value $\sigma_{c}$ we have in
practice pure fermionic star. On the figure it is seen that from
small values of $\mu_{c}$ to values near beyond the peak the rest
mass is greater than the total mass of the star which means that
the star is potentially stable.

\begin{figure}[ht]
\centerline{\psfig{figure=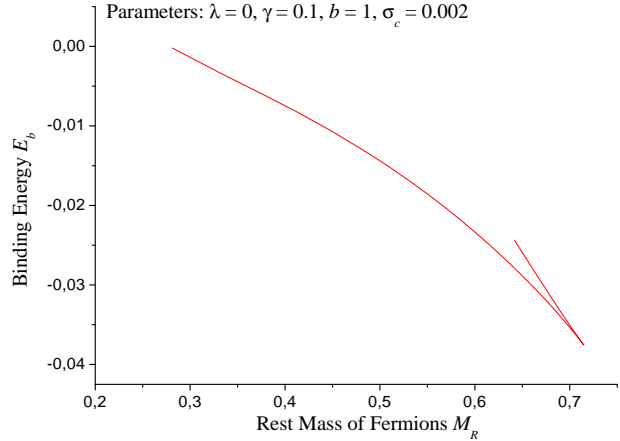,width=4.in}} \caption{The binding
energy $E$ versus the rest fermion mass $M_{RF}$} \label{bind}
\end{figure}

On Fig. \ref{bind} the binding energy of the star $E_{b}=M-M_{RB}-M_{RF}$ is drawn
against the rest mass of fermions $M_{RF}$ for $\lambda =0$, $\gamma =0.1$, $%
b=1$ and $\sigma _{c}=0.002$.
Fig. \ref{bind} is
actually a bifurcation diagram. With increasing the central value of function
$\mu(x)$
one meets a cusp. The appearance of a cusp shows that the stability of the
star changes - one perturbation mode develops instability. Beyond the cusp the
star is unstable and may collapse eventually forming black hole. The
corresponding physical results for pure boson stars are
considered in our recent paper \cite{FYBT2}.

\section*{Acknowledgments}
We are grateful to Prof. I.V. Puzynin (JINR, Dubna, Russia) for helpful discussion.

\end{document}